\documentclass[12pt,leqno]{article}
\usepackage{amsmath,amssymb,amscd,latexsym,bm}

\usepackage[all]{xy}

\newcommand{\C}{{\mathbb{C}}}

\newcommand{\Q}{{\mathbb{Q}}}
\newcommand{\oQ}{\overline{\Q}}
\newcommand{\R}{{\mathbb{R}}}

\newcommand{\X}{\mathbb{X}}
\newcommand{\Z}{{\mathbb{Z}}}

\newcommand{\Exp}{\mathrm{Exp}}

\newcommand{\id}{\mathrm{id}}

\newcommand{\spec}{\mathrm{spec}\,}

\newcommand{\Aut}{\mathrm{Aut}}

\newcommand{\End}{\mathrm{End}}

\newcommand{\Gal}{\mathrm{Gal}}
\newcommand{\GL}{\mathrm{GL}}

\newcommand{\Hom}{\mathrm{Hom}}

\newcommand{\emm}{{\mathfrak{m}}}

\newcommand{\eo}{\mathfrak{o}}

\newcommand{\onF}{\overline{F}}

\newcounter{aufz}
\newcommand{\beweisende}{\hspace*{\fill}$\Box$}

\newcommand{\verk}{\mbox{\scriptsize $\,\circ\,$}}

\begin{document}
\title{A remark on the vanishing of Higgs fields in the $p$-adic Simpson correspondence}
\author{Christopher Deninger\footnote{Funded by the Deutsche Forschungsgemeinschaft (DFG, German Research Foundation) under Germany's Excellence Strategy EXC 2044--390685587, Mathematics M\"unster: Dynamics--Geometry--Structure and the CRC 1442 Geometry: Deformations and Rigidity} \and Deepak Kamlesh\footnote{Supported by MPIM-Bonn-2025.}}
\date{}
\maketitle
\begin{abstract}
We give a condition on a $p$-adic representation of the fundamental group of a curve over $\overline{\Q}_p$ which ensures that under the $p$-adic Simpson correspondence the Higgs field vanishes. 
\end{abstract}
\section{Introduction} \label{sec:1}
The classical non-abelian Hodge correspondence on a compact Riemann surface gives an equivalence of categories between irreducible complex representations of the fundamental group and stable Higgs bundles of degree zero. The Higgs field is zero if and only if the representation is equivalent to a unitary representation.

The $p$-adic analogue of the non-abelian Hodge correspondence began with the pioneering work of Faltings \cite{F} and (for zero Higgs field) with \cite{DW1}. Since then the field has advanced significantly but the theory is not as complete as in the complex case. For example, according to the main result of \cite{X2} the Higgs bundles corresponding to $p$-adic representations of the fundamental group can be characterized in terms of their reduction behaviour (potentially strongly semistable). Such bundles are semistable of slope zero and for rank $\ge 2$ it is an important open question whether the converse is true. For a while this question seemed to have been settled in the negative by an ingenious counterexample but there was a subtle problem with the argument. The other basic open problem is to find a $p$-adic analogue of the unitarity condition on the representation of the fundamental group that corresponds to a vanishing Higgs field. The present note is a small contribution to this question in the following arithmetic case.

We fix an algebraic closure $\oQ_p$ of $\Q_p$ and set $\C_p = \widehat{\oQ}_p$. Let $X / \oQ_p$ be a smooth projective (connected) curve  of genus $g \ge 2$. The canonical embedding $\oQ_p \to B^+_{dR} / \xi^2$ determines a $B^+_{dR} / \xi^2$-lift $\X$ of $X$ which is one parameter on which the $p$-adic Simpson correspondence on $X_{\C_p} = X \otimes \C_p$ depends in the improved version of \cite{H}. The other parameter is the choice of an exponential $\Exp$ which is a continuous homomorphism splitting the logarithm $\log : 1 + \emm_{\C_p} \twoheadrightarrow \C_p$. For a point $x \in X (\oQ_p)$ the $p$-adic Simpson correspondence of \cite{H} restricts to a fully faithful functor
\[
S_{\Exp} : \begin{Bmatrix}
\text{continuous representations} \\ \rho : \pi_1 (X,x) \to \GL_r (\C_p) 
\end{Bmatrix} \longrightarrow
\begin{Bmatrix}
\text{Higgs bundles $( E, \theta)$ on $X_{\C_p}$} \\
\text{of rank $r$ with $E$ semistable}\\
\text{of slope zero} \end{Bmatrix} \; .
\]
Let $(X_K , x_K \in X_K (K))$ be a model of $(X ,x \in X (\oQ))$ over a finite extension $K / \Q_p$ with $K \subset \oQ_p$. Thus $X_K$ is a curve over $K$ with an isomorphism $X = X_K \times_{\spec K} \spec \oQ_p$ over $\oQ_p$, and under the inclusion $X_K (K) \subset X_K (\oQ_p) = X (\oQ_p)$ the point $x_K$ corresponds to $x$. The Galois group $G_K = \Gal (\oQ_p/ K)$ acts on $X$ from the left if $\sigma \in G_K$ corresponds to the automorphism of schemes $\id_{X_K} \times_{\spec K} \spec (\sigma^{-1})$. This induces a $G_K$-action on $X_K (\oQ_p) = X (\oQ_p)$ which keeps $x$ fixed. For $\sigma \in G_K$ we get an induced automorphism $\sigma_* = (\id_{X_K} \times_{\spec K} \spec (\sigma^{-1}))_*$ of $\pi_1 (X,x)$ and this leads to a left $G_K$-action on $\pi_1 (X,x)$. On the other hand $\sigma \in G_K$ defines a continuous group automorphism of $\GL_r (\C_p)$. For any continuous representation $\rho : \pi_1 (X,x) \to \GL_r (\C_p)$ we thus obtain a representation $^{\sigma} \rho$ as the composition
\[
^{\sigma} \rho : \pi_1 (X,x) \xrightarrow{\sigma^{-1}_*} \pi_1 (X,x) \xrightarrow{\rho} \GL_r (\C_p) \xrightarrow{\sigma} \GL_r (\C_p) \; .
\]
In this way the choice of the model $(X_K , x_K)$ gives a left action of $G_K$ on the set of representations $\rho$.  For another model over another finite extension $K' / \Q_p$ the two actions by $G_K$ and $G_{K'}$ on $\pi_1 (X,x)$ and its continuous representations coincide on an open subgroup of $G_K \cap G_{K'}$ in $G_{\Q_p}$. This holds since the models become isomorphic over a finite extension of $KK'$ in $\oQ_p$. Thus condition 1) of the theorem below is well defined. A stable vector bundle on a smooth projective curve over an algebraically closed field of characteristic zero is called \'etale-stable if its pullback to any finite \'etale Galois cover is stable. By \cite[Theorem 4.9]{W}, \'etale-stable bundles of rank $r \ge 2$ are generic. More precisely their locus in the moduli space of stable bundles is open with complement of codimension $\ge 2$. The notion of smallness of a representation $\rho$ is recalled in section \ref{sec:2}. 

\textbf{Theorem} \textit{For a continuous representation $\rho : \pi_1 (X,x) \to \GL_r (\C_p)$ let $(E, \theta) = S_{\Exp} (\rho)$ be the corresponding Higgs bundle. Assume either (i) $\rho$ is small and $E$ is stable or (ii) $E$ is \'etale-stable. Then the following assertions are equivalent. \\
1) For all $\sigma$ in some open subgroup of $G_{\Q_p}$ the representations $^{\sigma} \rho = \sigma \verk \rho \verk \sigma_*^{-1}$ and $\rho$ of $\pi_1 (X,x)$ are equivalent.\\
2) $\theta = 0$ and $E$ has a model over $X$.\\
In this case the matrices conjugating $\rho$ into $\rho^{\sigma}$ can be chosen to depend continuously on $\sigma$.}

If the vector bundle $E$ on $X_{\C_p}$ is stable then the representation $\rho$ is irreducible. The implication 2) $\Rightarrow$ 1) holds without assuming that $E$ is stable. This follows by combining \cite[Theorem 36]{DW1} with \cite{X1}. The Galois action is also studied in \cite[Theorem 25]{DW1}. For characters $\rho : \pi_1 (X,x) \to \C^{\times}_p$ a full characterization is known for when the Higgs field on the corresponding line bundle is zero. This happens if and only if $\rho$ can be approximated pointwise (or equivalently uniformly) by a sequence of characters $\rho_n$ which satisfy condition 1) i.e. for which there are open subgroups $G_n$ with $^{\sigma} \rho_n = \rho_n$, see \cite[Theorem 19]{DW2} if $X$ has good reduction and use \cite[Theorem 4.1]{H2} in general. 

A continuity argument for the $p$-adic Simpson correspondence using \cite[Theorem 1.1.1]{HX} shows that if a sequence of continuous representations $\rho_n : \pi_1 (X,x) \mapsto \GL_r (\C_p)$ as in 1) of the Theorem converges pointwise (or equivalently uniformly) to a continuous representation $\rho : \pi_1 (X,x) \to \GL_r (\C_p)$ then the Higgs bundle $S_{\Exp} (\rho)$ has vanishing Higgs field. It seems possible that as for rank one the existence of such a sequence $\rho_n \to \rho$ may also be necessary for $\theta = 0$ under the assumptions on $E$ in the theorem.

There does not seem to be an archimedian analogue of the theorem. For characters $\rho : \pi_1 (X (\C) , x) \to \C^{\times}$ on a curve $X / \R$ with $x \in X (\R)$, Galois equivariance means $\onF_{\infty *} (\rho) = \rho$. Here $F_{\infty} : X (\C) \to X (\C)$ is complex conjugation and $\onF_{\infty *} (\rho)$ is the composed representation, where $c = -$ denotes complex conjugation
\[
\onF_{\infty*} : \pi_1 (X,x) \xrightarrow{F_{\infty *}} \pi_1 (X,x) \xrightarrow{\rho} \C^{\times} \xrightarrow{c} \C^{\times} \; .
\]
Since $\onF_{\infty}$ corresponds to $F_{dR} = \id \otimes c$ on $X \otimes_{\R} \C$, Galois equivariant characters correspond to Higgs bundles $(L , \theta)$ where both the line bundle $L$ and the $1$-form $\theta$ are defined over the real algebraic curve $X$. This is easily checked using the explicit description of the Simpson correspondence for characters. There is thus no reason for $\theta$ to vanish.  

\textbf{Acknowledgement} For small representations, the authors discussed the main idea in this note long ago. The first author would like to thank Ben Heuer warmly for a very helpful e-mail exchange and several interesting comments. Thanks are also due to the referee for suggestions to improve the exposition. 

\section{Proof of the Theorem} \label{sec:2}

As explained above the choice of a model $(X_K , x_K)$ for $(X,x)$ over a finite extension $\Q_p \subset K \subset \oQ_p$ determines an action by $\sigma \in G_K$ on $X = X_K \times_{\spec K} \spec \oQ_p$ via $\id_{X_K} \times_{\spec K} \spec (\sigma^{-1})$. By the same formula, the inclusions $K \subset \oQ_p \subset B^+_{dR} / \xi^2$ determine a $G_K$-action on the canonical lift
\[
\X := X \times_{\spec \oQ_p} \spec (B^+_{dR} / \xi^2) = X_K \times_{\spec K} \spec (B^+_{dR} / \xi^2) \; .
\]
The natural morphism $\X \to X$ is then $G_K$-equivariant. The other parameter for the $p$-adic Simpson correspondence in \cite{H} is the choice of a continuous splitting $\Exp$ of the $p$-adic logarithm $\log : 1 + \emm_{\C_p} \to \C_p$. While $\log$ is $G_{\Q_p}$-equivariant, it is known that no splitting $\Exp$ can be equivariant for any open subgroup of $G_{\Q_p}$ and hence of $G_K$. This makes the relation between the Higgs bundles corresponding to a representation $\rho : \pi_1 (X,x) \to \GL_r (\C_p)$ and its conjugate representation $^{\sigma}\rho$ unclear. However this problem disappears if we assume that $\rho$ is small. Every continuous representation of a profinite group on $\C^r_p$ leaves an $\eo_{\C_p}$-lattice $\Gamma \subset \C^r_p$ invariant and small means that $\rho$ takes values in $1 + p^{\beta} \End_{\eo_{\C_p}} (\Gamma) \subset \Aut_{\eo_{\C_p}} (\Gamma)$ for some $\beta > 2 / (p-1)$ and some lattice $\Gamma$. In this case the Higgs bundle $(E, \theta)$ corresponding to $\rho$ is independent of the choice of $\Exp$: In the $p$-adic Simpson correspondence, $\Exp$ is used to construct exponential maps on certain Lie algebras which extend the usual exponential maps given by the standard convergent series on small enough subspaces. For small representations $\rho$ the relevant quantities already lie in the domain of the usual exponential maps and no uncanonical extension via a choice of $\Exp$ is necessary. Historically, Faltings dealt with the small case first and introduced a splitting $\Exp$ only later to deal with the general case. Since the usual exponential maps are Galois-equivariant it turns out that for a small representation $\rho$ and $\sigma \in G_K$ the Higgs bundle $^{\sigma} (E , \theta)$ corresponding to $^{\sigma} \rho$ is given as follows. Let $\sigma : X_{\C_p} \to X_{\C_p}$ be the automorphism $\id_{X_K} \times \spec (\sigma^{-1})$. Then
\[
^{\sigma} E = E \times_{X_{\C_p} , \sigma} X_{\C_p} \; .
\]
We have $\Omega^1_{X_{\C_p} / \C_p} = \Omega^1_{X_K / K} \otimes_K \C_p$ and hence conjugating $\theta : E \to E \otimes \Omega^1_{X_{\C_p} / \C_p}$ by $\sigma$, we obtain another Higgs field on $X_{\C_p}$
\[
^{\sigma} \theta : \,^{\sigma}E \longrightarrow \, ^{\sigma} E \otimes \Omega^1_{X_{\C_p} / \C_p} \; . 
\]
This not yet the Higgs field for $^{\sigma} \rho$. There is a twist. Let $\chi : G_{\Q_p} \to \Z^{\times}_p$ be the cyclotomic character, then we have
\[
^{\sigma} (E , \theta) = (\,^{\sigma} E , \, ^{\sigma} \theta \, \chi (\sigma)^{-1}) \; .
\]
One usually indicates this behaviour under Galois by writing $\theta \in \Hom (E , E \otimes \Omega^1_{X_{\C_p} / \C_p}) (-1)$. 

After these preliminaries, we can now show the remaining implication 1) $\Rightarrow$ 2). If $^{\sigma} \rho$ is equivalent to $\rho$ for all $\sigma$ in an open subgroup $H$ of $G_K$ then the Higgs bundle $(\,^{\sigma} E , \, ^{\sigma} \theta)$ and $(E, \theta)$ are isomorphic i.e. there is an isomorphism $\psi_{\sigma} : \,^{\sigma} E \to E$ over $X$ such that the following diagram commutes
\begin{equation}
\label{eq:1}
\xymatrix{
^{\sigma} E \ar[r]^-{^{\sigma}\theta} \ar[d]^{\wr}_{\psi_{\sigma}} & ^{\sigma} E \otimes \Omega^1_{X_{\C_p} / \C_p} \ar[d]^{\psi_{\sigma} \otimes \id} \\
E \ar[r]^-{\chi (\sigma) \theta} & E \otimes \Omega^1_{X_{\C_p} / \C_p}
}
\end{equation}
Such diagrams appeared earlier as ``enhanced Higgs bundles'' in the works \cite{MW}, \cite{T}, \cite{He} as Ben Heuer pointed out. By a theorem of Tate the fixed field $L = \C^H_p$ is finite over $\Q_p$. Let $M^0_K$ be the coarse moduli space of stable bundles of slope zero on the curve $X_K$. It is a quasiprojective variety over $K$. Since $^{\sigma} E$ and $E$ are isomorphic for $\sigma \in H$ we have equalities of isomorphism classes $^{\sigma} [E] = [\,^{\sigma} E] = [E]$ and hence
\[
[E] \in M^0_K (\C_p)^H = M^0_K (\C^H_p) = M^0_K (L) \subset M^0_K (\oQ_p) \; .
\]
It follows from the defining property of a coarse moduli space that $E \cong E_0 \otimes_{\oQ_p} \C_p$ where $E_0$ is a stable bundle on $X$. Note that we do not claim that $E_0$ is defined over $X_K \otimes_K L$.

After choosing a model $E_N$ of $E_0$ over some finite extension $N \supset K$ and replacing $H$ by an open subgroup of $H$ which fixes $N$, we have natural (transitive) identifications $^{\sigma} E = E$ for all $\sigma \in H$ over $^{\sigma} X = X$. Diagram \eqref{eq:1} now becomes the diagram
\begin{equation}
\label{eq:2}
\xymatrix{
E \ar[r]^-{^{\sigma}\theta} \ar[d]_{\psi_{\sigma}} & E \otimes \Omega^1_{X_{\C_p} / \C_p} \ar[d]^{\psi_{\sigma} \otimes \id} \\
E \ar[r]^-{\chi (\sigma) \theta} & E \otimes \Omega^1_{X_{\C_p} / \C_p} \; .
}
\end{equation}
Since we assume that $E$ is stable, we know that the isomorphism $\psi_{\sigma}$ is multiplication by a scalar. Hence \eqref{eq:2} implies that $^{\sigma} \theta = \chi (\sigma) \theta$ for all $\sigma \in H$. Writing $X_N = X_K \otimes_K N$ it follows that $\theta$ is an element of
\[
(\Hom_{X_N} (E_N , E_N \otimes \Omega^1_{X_N / N}) \otimes_N \C_p (-1))^H = \Hom_{X_N} (E_N , E_N \otimes \Omega^1_{X_N / N}) \otimes_N \C_p (-1)^H \; .
\]
By Tate's theorem, $\C_p (-1)^H = 0$ and hence $\theta = 0$. 

Now consider the general case, where we do not assume that the representation $\rho$ is small. By a well known argument, we can make $\rho$ small: Consider the exact sequence
\[
1 \longrightarrow 1 + p^{\beta} M_r (\eo_{\C_p}) \longrightarrow \GL_r (\eo_{\C_p}) \xrightarrow{\,\pi \,} \GL_r (\eo_{\C_p} / p^{\beta} \eo_{\C_p}) \longrightarrow 1 \; .
\]
Since $\GL_r (\eo_{\C_p} / p^{\beta} \eo_{\C_p})$ carries the discrete topology, the image of $\pi \verk \rho$ is finite. The kernel is therefore an open normal subgroup of $\pi_1 (X,x)$ which corresponds to a finite \'etale Galois cover $f : Y \to X$ of $X$ by a smooth projective curve $Y / \oQ_p$. By construction, the restriction $\rho' = \rho \, |_{\pi_1 (Y,y)}$ is a small representation of $\pi_1 (Y,y)$ for any point $y \in Y (\oQ_p)$ over $x$. Since $(Y,y)$ and $f$ are defined over a finite extension of $\Q_p$, it follows that for small enough open subgroups $H \subset G_{\Q_p}$ both $(Y,y)$ and $f$ are $H$-invariant. Hence $\pi_1 (Y,y)$ is an $H$-invariant subgroup of $\pi_1 (X,x)$ for small $H$ and therefore $^{\sigma} \rho'$ is equivalent to $\rho'$ for all $\sigma$ in a neighborhood of the identity in $G_{\Q_p}$. In general the Higgs bundle $(E' , \theta') = S_{\Exp} (\rho')$ is a twisted inverse image $f^{\circ} (E , \theta)$. In fact, in our arithmetic case we have $(E' , \theta') = f^* (E , \theta) = (f^* (E) , f^* (\theta))$ the ordinary inverse image. This holds because $X,Y$ and $f$ have \textit{compatible} lifts to $B^+_{dR} / \xi^2$ via $\oQ_p \to B^+_{dR} / \xi^2$. In order to show that $\theta = 0$ is suffices to prove that the Higgs field $\theta' = f^* (\theta)$ on $f^* (E)$ vanishes. This follows from the result in the small case since $f^* (E)$ is stable because we assumed that $E$ was \'etale-stable. \beweisende

\end{document}